\documentclass[12pt,reqno]{amsart}
\usepackage{amssymb}
\usepackage{epsf}
\usepackage{amsxtra}

\theoremstyle{plain}
\newtheorem{Thm}{Theorem}

\newtheorem{Lem}{Lemma}

\theoremstyle{definition}

\theoremstyle{remark}

\newtheorem*{Ex}{Example}

\numberwithin{equation}{section}

\begin{document}

\title{Peg Jumping for Fun and Profit}

\date{\today}

\author{David~M. Bradley}
\address{Department of Mathematics \& Statistics\\
         University of Maine\\
         5752 Neville Hall
         Orono, Maine 04469-5752\\
         U.S.A.}
\email[David Bradley]{bradley@math.umaine.edu,
dbradley@member.ams.org}
\author{Hugh Thomas}
\thanks{The second author was
partially supported by a grant from the Natural Science and Engineering
Research Council (Canada)}
\address{Department of Mathematics \& Statistics\\
University of New Brunswick\\
Fredericton, NB\\
E3B 5A3 Canada}
\email[Hugh Thomas]{hugh@math.unb.ca}



\begin{abstract} We consider the problem of determining the
minimum number of moves needed to solve a certain one-dimensional
peg puzzle.  Let $N$ be a positive integer. The puzzle apparatus
consists of a block with a single row of $2N+1$ equally spaced
holes which, apart from the central hole, are occupied by an equal
number $N$ of red and blue pegs.  The object of the puzzle is to
exchange the colors of the pegs by a succession of allowable
moves.  Allowable moves are of two types: a peg can be shifted
from the hole it occupies into the empty hole adjacent to it, or a
peg can jump over an adjacent peg into the empty hole.  We exhibit
a sequence of $N^2+2N$ moves that solves the puzzle, and prove
that no solution can employ fewer moves.
\end{abstract}

\maketitle

\section{Introduction}\label{sect:Intro}
We are going to begin by describing a puzzle.
Let $N$ be a positive integer. Consider the one-dimensional peg
puzzle whose apparatus consists $N$ red pegs, $N$ blue pegs, and a
block with a single row of $2N+1$ equally spaced holes. Initially,
the $N$ red pegs occupy the rightmost $N$ holes, the $N$ blue pegs
occupy the leftmost $N$ holes, and the center hole is empty. It is
required to exchange the red pegs and the blue pegs so that the
red pegs occupy the leftmost $N$ holes, and the blue pegs occupy
the rightmost $N$ holes.  We would also like to carry out this
task using as few moves as possible. Allowable moves are of two
types: a \emph{step}, which consists of one peg moving from an
occupied hole into the adjacent (necessarily unique) unoccupied
hole; and a \emph{jump}, which consists of one peg jumping over a
single adjacent peg into the empty hole.

The puzzle as just described with $N=5$ was marketed by
International Games of Canada Ltd., Mississauga, Ontario, under
the name of ``Brainbuster,'' back in the 1970s.
In this paper, we provide a sequence of $N^2+2N$ moves which solves the
puzzle for general $N$, and we show that no shorter sequence of moves
is sufficient.

Before continuing on to read our results,
the reader may prefer to attempt some small cases.
$N=2$ is pretty easy, and $N=3$ is difficult enough to give the
flavour of the general solution.

%
%

\section{A Solution Consisting of $N^2+2N$ Moves}\label{Sect:UpperBound}

To describe our solution requires some notation.  Encode a step to
the right made by a blue peg by $S$, and a step to the left made
by a red peg by $s$.  Similarly, encode a jump to the right made
by a blue peg by $J$, and a jump to the left made by a red peg by
$j$.  The solution takes a slightly different form depending on
the parity of $N$.
\begin{Thm}\label{thm:solnseq}
Let $n$ be a positive integer.  If $N=2n$, then
\begin{equation}
\label{even}
   \bigg(\prod_{k=1}^n Sj^{2k-1}sJ^{2k}\bigg)
   \bigg(\prod_{k=1}^n sj^{2n-2k+1}SJ^{2n-2k}\bigg)
\end{equation}
is a sequence of $N^2+2N$ moves that solves the puzzle. If
$N=2n-1$, then
\begin{equation}
\label{odd}
   \bigg(\prod_{k=1}^{n-1}Sj^{2k-1}sJ^{2k}\bigg)Sj^{2n-1}S
   \bigg(\prod_{k=1}^{n-1}J^{2n-2k}sj^{2n-2k-1}S\bigg)
\end{equation}
is likewise a sequence of $N^2+2N$ moves that solves the puzzle.
\end{Thm}

\begin{Ex} For $N=1,2,3,4,5$ we obtain the respective solution
sequences $SjS$, $SjsJJsjS$, $SjsJJSjjjSJJsjS$,
$SjsJJSjjjsJJJJsjjjSJJsjS$ and
\[
   SjsJJSjjjsJJJJSjjjjjSJJJJsjjjSJJsjS.
\]
\end{Ex}

The solutions defined in the statement of  Theorem~\ref{thm:solnseq}
can be reformulated
in a manner which emphasizes their symmetry.

If $N=2n$, then ~\ref{even} can be refomulated as:
\begin{equation}\label{eqeven}
\bigg(\prod_{k=1}^{n-1}Sj^{2k-1}sJ^{2k}\bigg)
   Sj^{2n-1}sJ^{2n}sj^{2n-1}S\prod_{k=n-1}^1
   J^{2k}sj^{2k-1}S.
\end{equation}
If $N=2n-1$, then ~\ref{odd} can be reformulated as:
\begin{equation}\label{eqodd}
\bigg(\prod_{k=1}^{n-1}Sj^{2k-1}sJ^{2k}\bigg)
   Sj^{2n-1}S\bigg(\prod_{k=n-1}^1 J^{2k}sj^{2k-1}S\bigg).
\end{equation}

\noindent
{\bf Proof of Theorem~\ref{thm:solnseq}.}
We consider the case where $N=2n-1$ is odd, leaving the
case where $N$ is even to the reader.  We provide
illustrations for the case $N=5$.
We start
with the position $B^NOR^N$.
$$\epsfbox{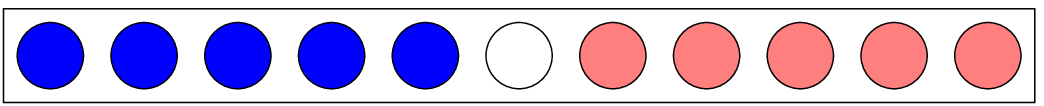}$$
We apply $Sj$, obtaining
$B^{N-1}RBOR^{N-1}$.
$$\epsfbox{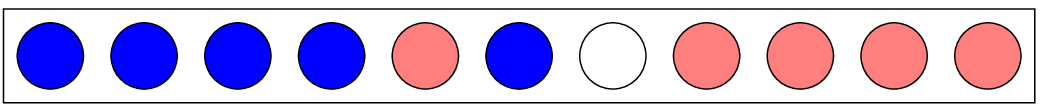}$$
Then we apply $sJJ$, obtaining
$B^{N-2}0RBRBR^{N-2}=B^{N-2}O(RB)^2R^{N-2}$.
$$\epsfbox{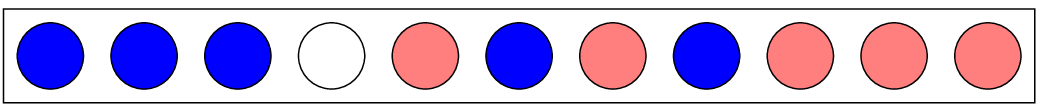}$$
We have now applied the first term in the product
(\ref{eqodd}).
It is easily checked that after the application of $m$ terms of the product,
we obtain
$B^{N-2m}O(RB)^{2m}R^{N-2m}$.  After all $n-1$ terms of the product, we obtain
$BO(RB)^{2n-2}R$:
$$\epsfbox{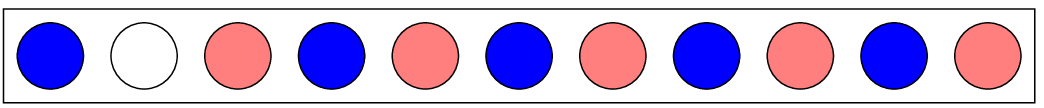}$$
We then apply $Sj^{2n-1}S$, obtaining $R(BR)^{2n-2}OB$:
$$\epsfbox{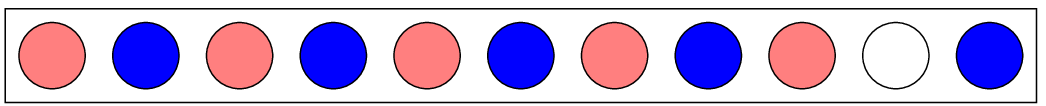}$$
Observe that this is the reversal of the previous result, the result of
applying only the first
product.  The second product then acts in the opposite way
to the first product, and the final result is
$R^NOB^N$, as desired:
$$\epsfbox{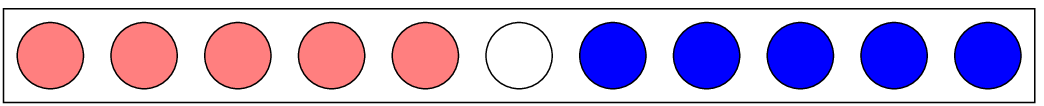}$$
This completes the proof that the sequence of moves given in the statement
of the theorem is indeed a solution.

Now that we have solutions to the puzzle, we can count up the number of
moves which they require.  Using whichever of
(\ref{eqeven}) or (\ref{eqodd}) applies,
it is easy to see that $2N$ steps are used in either case, while
the number of jumps is
$$\left(\sum_{i=1}^{N-1} i \right) + N +\left(\sum_{i=1}^{N-1} i\right)
= \frac{N(N-1)}{2} + N + \frac{N(N-1)}2 = N^2.
$$
Thus the solutions given in the statement of the theorem
use a total of $N^2+2N$ moves, as desired.

\section{$N^2+2N$ moves are necessary}

This section is devoted to the proof of the following theorem:

\begin{Thm}\label{thm:nec} Any solution to the peg-jumping problem requires at least
$N^2+2N$ moves.
\end{Thm}

\noindent
{\bf Proof of Theorem~\ref{thm:nec}.}
Let us refer to the {weight} of a position as the total distance to the
right that the blue pegs have moved plus the total distance to the left
that the red pegs have moved.  The weight of the initial position is zero,
and the weight of the final position is $2N(N+1)$.

The weight increases by 2 if a blue peg jumps to the right, or if a
red peg jumps to the left.  It increases by 1 if a blue peg steps to the
right or a red peg steps to the left.  The other possible moves
decrease the weight.

Since on any move, the weight increases by at most 2, clearly
$$\frac{2N(N+1)}2=N^2+N$$
moves are necessary.  However, this isn't quite good enough for our purposes.
This bound would only be attained if all the pieces only ever jumped, and
no such solution is possible; in fact, it is easy to verify that, from
any position, no more than $N$ consecutive jumps which increase the weight
are possible.

We say that a pair of pegs {\it crosses} when one jumps over the other.
If we pick a red peg and a blue peg, since the red peg starts to the
right of the blue peg and ends to its left, this pair of pegs must cross
an odd number of times; in particular, this pair of pegs must cross at
least once.  Let us refer to a move as a {\it first cross} if it is a move
on which a red peg and a blue peg cross for the first time.
In any solution, there
must be such a cross for each such pair;
thus, there must be $N^2$ first crosses.  These $N^2$ first crosses have
the total effect of increasing the weight of the position by $2N^2$, leaving
a weight increase of $2N$ which must be coming from other moves.

We claim that the average increase in weight over all the moves other than
the first crosses is at most 1.  This is sufficient to prove the theorem,
because it means that the remaining weight increase of $2N$ will require
at least $2N$ moves in addition to the first crosses to accomplish this
weight gain, which will establish the theorem.

So now we must prove the following lemma:
\begin{Lem}\label{lem:avg} In any solution to the puzzle,
the average weight gain over all moves
other than first crosses is at most 1.
\end{Lem}

\noindent
{\bf Proof.} Our procedure will be to group the moves other than the first
crosses into small groups, each of which groups we show to have average
weight gain of at most 1.

If a red peg $R$ jumps a blue peg $B$ other than at their first crossing,
group together all the crossings of these two pegs.  We have already
argued that after removing the first crossing there will be an even number,
half going in one direction, half in the other, so the total weight gain
for all these crossings will be zero.

A move consisting of stepping a red piece to the left by 1 can go in a group
by itself; it adds 1 to the weight, so it has average weight 1.  We proceed
similarly for stepping a blue piece to the right.

What moves which add to the weight remain?  Only jumping a red piece over
another red piece, moving left, or jumping a blue piece over a blue piece,
moving to the right.  Either of these moves adds 2 to the weight, so we
must pair them with some other moves.  We will shortly describe which moves
we want to use.

We pause this analysis to consider how the game looks from
the point of view of a single peg $P$.  We consider that a peg $P$
only notices when it is being moved or jumped, or when someone jumps over it.
So on a given turn, the peg $P$ will notice one of five things:
\begin{itemize}
\item[(A)] Someone jumped over $P$ proceeding left.
\item[(B)] $P$ moved or stepped left.
\item[(C)] Someone jumped over $P$ proceeding right.
\item[(D)] $P$ moved or stepped right.
\item[(E)] Nothing.
\end{itemize}

We now claim that, among turns where $P$ notices something happening,
turns of type (A) or (B) alternate with turns of type (C) or (D).  In
other words, if, on one turn, either (A) or (B) happened, the next time
$P$ notices something happen, that event will be either (C) or (D).
The reason for this is that after either (A) or (B) has happened, the
empty hole is to the right of $P$.  This will only change if (C) or
(D) happens, and while the empty hole remains to the right of $P$,
(A) and (B) are impossible.  A
similar argument applies if we start with a move of type (C) or (D).

We now consider how to pair up moves involving same-colour jumping.
We begin by describing a procedure which doesn't quite work, and then
we describe how to fix it so it really does work.  Suppose that
a red peg $S$ jumps a red peg $R$, proceeding to the left.  From the
point of view of $R$, this is a move of type (A).  Thus, the next
move $R$ notices will be a move of type (C) or (D).  Suppose it is
of type (C): $R$ notices someone jump to the right over him.  But this
piece jumping to the right over him must be $S$ (since the empty hole is
to the right of $R$, $S$ cannot have moved since jumping $R$).  We can
then pair up these two jumps, and the total weight gain of this pair is
zero.

Suppose now that the move after $S$ jumped over $R$ is of type (D).
On this move, $R$ moves rightwards, i.e. in a direction which decreases
weight.  Thus if we pair the jump of $S$ over $R$ with this backwards
move, the total weight increase will be at most 1 over the two moves, and
we are certainly satisfied.

The problem with this way of pairing up moves
is that there might be some same-colour jump which we want to pair up
with a following move of type (C) or (D), but that following move never
occurs because the puzzle is solved first.
Thus, we need to do something
slightly more complicated.

When $S$ jumps to the right over $R$, we take note of whether or not $R$ has
already had all the blue pegs cross it.  If not, then we pair the
jump with the next move which $R$ sees, as described above.  (As we already
saw, this move will consist of $R$ moving to the right
or else of $S$ jumping to
the right
over $R$, and since $R$ still has to see some blue peg cross it before the
puzzle can be solved, and this
can only happen {\it after} the next event $R$ sees, the puzzle cannot be
finished
before the event which we want to pair with the jump of $S$ over $R$.)

On the other hand, if $R$ has already seen all the blue pegs cross
it when $S$ jumps over it to the left,
we pair the jump with the
{\it previous} move which $R$ sees.  By the same argument as above,
 this move will be either
$S$ jumping to the right over $R$, or else $R$ moving to the right.  This
move occurs {\it after} the last blue peg jumped over $R$, so this move
is not already being paired with some other same-colour jump over $R$,
and thus we have paired up all the jumps of red pegs over red pegs.

We proceed similarly with jumps of blue pegs over blue pegs.
The only remaining
moves are ones which decrease the total weight, and we do not need to
group them.  This proves the lemma, and with it the theorem.
\end{document}